\documentclass[11pt]{amsart}

\usepackage{amsmath,amssymb,latexsym,amsthm,enumerate}
\usepackage{verbatim}
\usepackage{graphicx} 
\usepackage{hyperref} 
 
\usepackage{xcolor}

\def\NN{\mathbb{N}}

\def\RR{\mathbb{R}}

\def\F{\mathcal{F}}

\def\1{\mathbf{1}}
\def\0{\mathbf{0}}

\def\I{\mathbf{I}}

\newcommand{\te}{\widehat{e}}

\newtheorem{Thm}{Theorem}
\theoremstyle{definition}
\newtheorem*{Def}{Definition}

\newtheorem{As}{Assumption}

\def\BEN{\begin{enumerate}}  \def\BI{\begin{itemize}}
\def\EEN{\end{enumerate}}   \def\EI{\end{itemize}}
    \def\nn{\nonumber}

\def\mbb{\mathbb}  \def\mrm{\mathrm}
\def\mc{\mathcal}

\def\le{\left}
\def\ri{\right}
\def\te#1{\mathrm{e}^{#1}}   

\def\WH{\widehat}

  \def\d{\delta}   \def\th{\theta}

  \def\nn{\nonumber}   
     
\def\w{\omega} \def\q{\qquad} 

  \def\td{\text{\rm d}}
 
\numberwithin{equation}{section}

\newcommand{\exit}{{\mbox{\, \vspace{3mm}}}
\hfill\mbox{$\square$}}

\begin{document}

\title[Asymptotic independence of statistics of maximal increments]{Asymptotic 
independence of three statistics \\of maximal segmental scores}

\author{Aleksandar Mijatovi\'{c}}
\author{Martijn Pistorius}
\address{Department of Mathematics, Imperial College London}
\email{\{a.mijatovic,m.pistorius\}@imperial.ac.uk}

\keywords{Maximal segmental score, asymptotic independence, 
asymptotic overshoot, random walk, de-Poissonization.}

\thanks{We are grateful to the participants of the probability
and stochastic analysis seminars at TU Berlin, Oxford and Warwick
for comments which improved the paper. AM would like to thank Vlad 
Vysotsky for useful discussions.}

\subjclass[2000]{60G51, 60F05, 60G17}

\begin{abstract}
Let
$\xi_1,\xi_2,\ldots$
be an iid sequence with negative mean. The 
$(m,n)$-segment 
is the subsequence 
$\xi_{m+1},\ldots,\xi_n$
and its \textit{score} is given by
$\max\{\sum_{m+1}^n\xi_i,0\}$.
Let
$R_n$
be the largest score of any segment
ending at time 
$n$,
$R^*_n$
the largest score of any segment in the sequence
$\xi_{1},\ldots,\xi_n$,
and 
$O_x$
the overshoot of the score over a level
$x$
at the first epoch the score of such a size arises.
We show that, under the Cram\'er assumption 
on 
$\xi_1$,
asymptotic independence
of the statistics 
$R_n$,  $R_n^* -y$ and $O_{x+y}$ holds
as $\min\{n,y,x\}\to\infty$.  
Furthermore, we establish a novel Spitzer-type identity characterising 
the limit law $O_\infty$ in terms of the 
laws of $(1,n)$-scores.
As corollary we obtain: (1) a novel factorization of 
the exponential distribution as a convolution 
of 
$O_\infty$ 
and the stationary distribution of 
$R$;
(2)
if 
$y=\gamma^{-1}\log n$
(where $\gamma$ is the Cram\'er coefficient),
our results, together with the classical theorem of Iglehart~\cite{Iglehart}, 
yield the existence and explicit form of the joint weak limit 
of
$(R_n, R_n^* -y,O_{x+y})$.
\end{abstract}

\maketitle

\section{Introduction and the main result}
\label{sec:Intro}

Consider a sequence of iid random variables 
$\{\xi_i\}_{i\in\mathbb N}$
with negative mean
and denote by $S=\{S_n\}_{n\in\mathbb N^*}$ the
random walk corresponding to  
$\{\xi_i\}$:
$S_0\doteq0$ and $S_n\doteq \sum_{i=1}^n\xi_i$. 
For any $m<n$ with $n\in\NN$, $m\in\NN^*\doteq\NN\cup\{0\}$, 
the \textit{segmental score} of the $(m,n)$-segment
$\{\xi_i\}_{i=m+1}^n$ of $\{\xi_i\}$ 
is given by the maximum of the sum of the elements in the 
segment and zero
(as usual we denote $x^+\doteq \max\{x,0\}$, $x\in\mbb R$):
$$\left(\sum_{i=1+m}^n\xi_i\right)^+=(S_n-S_m)^+.$$
The notion of segmental scores arises naturally in several
areas of applied probability and statistics. For their application
in the study of DNA sequences see {\em e.g.}~\cite{Avery_Henderson} and~\cite{KarlinDembo}. 
Segmental scores also play an important role in sequential change point detection 
problems of mathematical statistics ({\em e.g.}~CUSUM test), see~\cite{Siegmund} and~\cite{Mous,Shiryaev}.
Moreover, in sequential analysis in the context of abortion epidemiology,
the maximal segmental score is proposed in~\cite{Levin_Kline} 
as a test statistic for the detection of a \textit{one-sided epidemic alternative} 
for the increase in the mean of a sequence of independent random variables
(see also~\cite{Commenges_etal} for related applications of 
the epidemic alternative in experimental neurophysiology).
For the role of the segmental scores in queueing theory see~\cite{AsmussenC,AsmussenAPQ,Iglehart}. 
It is of interest in all of these applications to quantify the fluctuations of the 
segmental scores. In the recent paper~\cite{MikoschRackauskas} 
(see also~\cite{MikoschMoser_13}) this problem is studied  under 
heavy-tailed step-size distributions, where an appropriate scaling 
of segmental scores is necessary for the analysis. In the case of 
an exponentialy thin positive tail, {\em i.e.}~under the Cram\'er assumption,
no scaling is required and the asymptotics of the fluctuations of the 
segmental scores can be analysed directly, which is the aim
of the present paper.



Two natural statistics measuring the fluctuations of the segmental scores 
are $R_n$, the largest score of any $(m,n)$-segment
({\em i.e.}~of any segment ending at time $n$), 
and $R_n^*$, the largest score of any of segment in $\{\xi_i\}_{i=1}^n$
(i.e. the largest score that has arisen up to time $n$).  More
precisely, for $n\in\mbb N$, we have 
\begin{eqnarray*}
R_n\, \doteq\, \max_{m\in\{0, \ldots, n-1\}}(S_{n} - S_m)^+ 
\qquad \text{and}\qquad
 R_n^*\, \doteq\,
\max_{m,k\in\{0, \ldots, n\}, m < k} (S_{k}-S_m)^+, 
\end{eqnarray*}
A third statistic quantifying the fluctuations of 
segmental scores is the first segmental score larger than
$x>0$, which is given by  $R_{H(x)}$ with $H(x)$ the first time 
an increment of the random walk larger than $x$ occurs:
\begin{equation*}
H(x) \,\doteq\,  \min\{k\in\mbb N: \text{$\exists\> m<k$
such that $S_{k}-S_m>x$}\}.
\end{equation*}
The main contribution of this paper is to give sufficient conditions for the three
statistics
$$
R_n, \qquad Q_{n,x} \,\doteq\, R^*_n - x, \qquad O_x \doteq R_{H(x)} - x,
$$
of the maximal increments of the walk 
$S$
to be asymptotically independent in the sense
that the joint CDF is asymptotically equal to the 
product of the marginal CDFs of the statistics:

\begin{Def}
A family of random vectors 
$\{(U_z^1, \ldots, U^d_z)\}_{z\in\mathcal Z}$
on a given probability space,
indexed by
$z\in\mathcal Z\subset[0,\infty)^l$,
$d,l\in\NN$,
is \textit{asymptotically independent} if
the joint CDF is asymptotically equal to a product 
of the CDFs of the components: i.e. 
for any
$a_i\in(-\infty,\infty]$, $i=1,\ldots,d$,  it holds
\begin{equation*}
P(U_z^1 \leq a_1, \ldots, U_z^d\leq a_d) = \prod_{i=1}^d P(U^i_z\leq a_i)
 + o(1) 
\qquad
\text{as $\min\{z_1,\ldots, z_l\}\to\infty$.}
\end{equation*}
\end{Def}
\noindent Our result states that the asymptotic independence of the three statistics 
above essentially holds under the Cram\'er assumption on the step-size distribution
(which in particular implies 
$E[\xi_1]<0$):
\begin{As}\label{A1}
The distribution of $\xi_1$ has finite mean, is non-lattice and
satisfies \emph{Cram\'{e}r's condition}, i.e. $E[\te{\gamma
\xi_1}]=1$ for some $\gamma\in(0,\infty)$, and 
$E[|\xi_1|\te{\gamma \xi_1}]$ is finite.
\end{As}

\begin{Thm}\label{thm2}
Under As.~\ref{A1}, 
the triplet $\{(R_n,Q_{n,y},O_{y+x})\}_{n\in\mbb N, x,y\in\mbb R_+}$ 
is asymptotically independent,
where 
$\RR_+\doteq[0,\infty)$.
Furthermore,
the following limit in distribution holds:
$O_x\stackrel{\mc D}{\longrightarrow} O_\infty$
as
$x\to\infty$,
where
$O_\infty$
is a non-negative distribution with the characteristic function
\begin{equation}\label{eq:oo}
E[\te{\mrm i\th\, O_\infty}] = \frac{\gamma}{\gamma - \mrm i
\theta} \cdot \exp\le\{\sum_{n=1}^\infty\frac{1}{n}\left(1 -
E\left[\te{\mrm i \th\, S_n^+}\right]\right)\ri\},
\qquad \text{for all}\>\>\th\in\mbb R.
\end{equation}
\end{Thm}
\noindent {\bf Remarks.} 
(i) A classical time-reversal argument 
implies that
$R_n$ 
and
$\max_{m\in\{0, \ldots, n\}}S_m$
have the same law 
for every $n\in\NN$.
Hence $R_n$
converges in distribution, as 
$n\uparrow\infty$,
to $S^*_\infty\doteq\sup_{n\in\NN}S_n$, 
which is finite (as $E[\xi_1]<0$ by As.~\ref{A1}) and follows a distribution characterised by Spitzer's identity 
(see~\cite[p.230]{AsmussenAPQ}). 

(ii) Note that Spitzer's identity~\cite[p.230]{AsmussenAPQ}
and a time-reversal argument imply that 
the second factor in~\eqref{eq:oo}
is equal to $1/E[\te{\mrm i\theta R_{\infty}}]$.
The asymptotic independence of Thm.~~\ref{thm2}
therefore yields the joint law of the weak limit 
$(R_\infty,O_\infty)$.
In particular, the limit law 
$R_\infty+O_\infty$
of the sum
$R_n+O_x$,
as
$\min\{x,n\}\to\infty$,
is characterised by the identity
$$E\left[\te{\mrm i\theta (R_{\infty}+O_\infty)}\right] =
\frac{\gamma}{\gamma - \mrm i \theta},\>
\>\>\forall\th\in\mbb R,\quad
\text{ and hence }\quad\gamma\cdot\le(R_\infty + O_\infty\ri)\sim\mathrm{Exp(1)}.$$
This establishes a novel factorization of the exponential distribution
$\mathrm{Exp(1)}$
as the convolution of the distribution of the asymptotic overshoot 
and the stationary distribution of a reflected random walk 
with step-size distribution satisfying As.~\ref{A1}.
Note further that, unlike in the Wiener-Hopf factorisation, here the supports of the 
factorising random variables are in general not disjoint.

(iii) Note that $Q_{n,x}$ does not admit
a non-degenerate weak limit along any sequence $(n,x)$ with $\min\{n,x\}\to\infty$.
A sufficient condition for the weak convergence 
of the statistic $Q_{n,x}$ is given in Iglehart~\cite[Thm.~2]{Iglehart}:
if $x(n) = \gamma^{-1} \log( K n)$,
for a certain positive constant 
$K$, then 
$\gamma\, Q_{n,x(n)}$ converges weakly 
to a Gumbel distribution as $n\uparrow\infty$.

(iv) The main technical fact established in this paper
is that asymptotically, as $\min\{x,y,n\}\to\infty$, the probability 
that $R$ crosses the level $x+y$ for the first time
during the excursion of $R$ away from 0 straddling time
$n$
vanishes
(see Section~\ref{ssec:asympi}).
This fact,
in conjunction with the independence of distinct excursions, 
essentially implies  the asymptotic independence in Thm.~\ref{thm2}.

\section{Proof of Thm.~\ref{thm2}}\label{sec:proof}
We start with the observation that $R$ is a reflected random walk 
and that $R^*$ and $H(x)$ may be represented as path-functionals of $R$:
\begin{equation*}
R_n= S_n - \min_{m\in\{0, \ldots, n\}} S_m,
\qquad 
 R_n^*= \max_{k\in\{1, \ldots, n\}} R_k
 \quad\text{and}\quad
H(x)=\min\{n\in\mbb N: R_n > x\}.
\end{equation*}
The proof of the asymptotic independence is based on (i) an application 
of asymptotic results for reflected L\'evy processes 
(reviewed in Section~\ref{subsec:Levy}) to the compound 
Poisson process obtained by subordinating $S$ by a Poisson 
process and (ii) asymptotic de-Poissonization results, 
which are established by splitting arguments (given in Section~\ref{ssec:asympi}). Finally, the proof of Eqn.~\eqref{eq:oo} is presented in Section~\ref{ssec:spit}.

\subsection{Maximal increments of L\'evy processes}
\label{subsec:Levy}

Let $X=\{X(t)\}_{t\ge 0}$
be a L\'{e}vy process under a probability 
measure
$P$, {\em i.e.}~a $P$-a.s. c\`adl\`ag process started at
$X(0)=0$
with stationary independent increments 
(refer to~\cite{Bertoin} for background on the theory of 
L\'evy processes). 
Denote by $Y=\{Y(t)\}_{t\ge 0}$ the L\'{e}vy process reflected at its running
infimum and by 
$\tau(x)$
the first time an increment of $X$
of size larger than 
$x\in\mbb R_+$
occurs:
\begin{equation*}
Y(t) \,\doteq\, X(t) - \inf_{0\leq s\leq t}X(s) = \sup_{0\leq s\leq t}\{X(t)-X(s)\},
\quad
\tau(x) \,\doteq\, \inf\{t\ge 0: Y(t) > x\} \text{ ($\inf\emptyset=\infty$}).
\end{equation*}
Let
$Y^*(t)\doteq\sup_{0\leq s\leq t}Y(s)$ 
be the supremum of
$Y$ up to time $t$.  
The three statistics of interest take the following
form for any 
$t, x \in \mbb \RR_+$:
\begin{equation}
\label{eq:Levy_Overshoot}
Y(t),\qquad 
M({t,x}) \doteq Y^*(t) - x,\qquad
Z(x)\doteq Y({\tau(x)}) - x.
\end{equation}

Denote by $\WH L$ a local time at zero of $Y$ and 
let $L$ be a local time  at zero
of the reflected process $\WH Y=\{\WH Y(t)\}_{t\ge 0}$ of
the dual 
$\WH X \doteq -X$,
$\WH Y(t) \doteq \sup_{0\leq s\leq t}X(s)- X(t).$
The ladder-time process
$L^{-1}~=~\{L^{-1}(t)\}_{t\ge 0}$
is equal to the right-continuous inverse of 
$L$ (see~\cite[Ch.~IV]{Bertoin} for a definition of local time and its inverse).
The ladder-height process
$H=\{H(t)\}_{t\ge 0}$ is given by $H(t) \doteq 
X(L^{-1}(t))$
for all $t\ge 0$ with $L^{-1}(t)$ finite and by
$H(t)\doteq+\infty$
otherwise.
Let $\phi$ be the Laplace
exponent of $H$,
\begin{equation}\label{eq:phi}
\phi(\theta)\doteq-\log E[\te{-\theta H(1)}\I_{\{H(1)<\infty\}}],
\qquad\text{for
any}\quad \theta\in\RR_+,
\end{equation}
where 
$\I_A$ denotes the indicator of a set $A$
and
$E[\cdot]$
is the expectation under
$P$.


\begin{As}\label{A2} 
The mean of $X(1)$ is finite,
the Cram\'{e}r condition 
for 
$X$
is satisfied, i.e.
\begin{equation}\label{cr}
\text{there exists a $\gamma\in(0,\infty)$ such that $E[\te{\gamma X(1)}]=1$,}
\end{equation}
$E[\te{\gamma X(1)}|X(1)|]<\infty$,
and either the 
L\'evy measure of
$X$
is non-lattice or
$0$
is regular for
$(0,\infty)$.
\end{As}
\noindent 
\textbf{Remark.} Display~\eqref{cr}
implies $E[X(1)]<0$, 
making
$Y$ 
(resp. $\WH Y$)
a recurrent 
(resp. transient) Markov process on 
$\RR_+$. Hence 
$\phi(0)>0$
and the stopping time $\tau(x)$ 
is a.s. finite for any 
$x\in\RR_+$,
so that  
$H$ is
a killed subordinator
and the overshoot
$Z(x)$ is a well-defined random variable.

Thm.~\ref{thm}, 
which provides a key
step in the proof of Thm.~~\ref{thm2}, 
was established in~\cite[Thm.~1,~Lem.~3]{MiPi_Levy}.

\begin{Thm}\label{thm}
Under As.~\ref{A2} the following statements hold:
\begin{enumerate}[(i)]
\item 
\label{enum:i}
The triplet
$\{(Y(t), Z(x+y), M({t,x}))\}_{t,x,y\in\mbb R_+}$
is asymptotically independent.
\item
\label{enum:ii}
The limit in distribution holds:
$Z(x) \stackrel{\mc D}{\longrightarrow} Z(\infty)$
as $x\to\infty$,
where 
\begin{equation}\label{eq:zinf}
E[\te{-v Z({\infty})}] = \frac{\gamma}{\gamma + v} \cdot \frac{\phi(v)}{\phi(0)}
\qquad\text{for all}\quad v\in\RR_+.
\end{equation}
\item 
\label{enum:iii}
It holds that
$\lim_{x\wedge t\to\infty}P(\WH L(t) = \WH L(\tau(x))) =0$ 
and
for any $\delta_1,\delta_2\in[0,1/4)$ we have 
\begin{eqnarray}
\label{eq:Lemma_iv_ineq}
\limsup_{x\wedge t\to\infty}
P\le(\WH L(t(1-\delta_1)) \leq \WH L(\tau(x))\leq \WH L(t(1+\delta_2))\ri)\leq
\frac{8}{\te{}}\max\{\delta_1,\delta_2\}.
\end{eqnarray}
\end{enumerate}
\end{Thm}

\subsection{Proof of the asymptotic independence}\label{ssec:asympi} 
Let the random walk
$\{S_n\}_{n\in\NN^*}$ 
and an independent Poisson process
$\{N(t)\}_{t\geq0}$
with unit intensity rate
be defined on the same probability space, 
and define a compound Poisson process
$X=\{X(t)\}_{t\geq0}$ by
\begin{equation}
\label{eq:Emebdding}
X(t)\doteq S_{N(t)},\quad t\geq0.
\end{equation}
For any
$t >0$,
let $[t]\doteq\max\{n\in\NN^*:n<t\}$
denote the largest integer which is smaller than
$t$, and set $[0]\doteq 0$. 
For any $t, x,y>0$ let $A,B,C$ and $A',B',C'$ be the sets
\begin{eqnarray}
\nn && A \doteq \{Y(t)\leq w\}, \qquad B \doteq  \{Z(x+y)\leq v\}, \qquad C \doteq  \{M(t,y)> z\},\\
\label{abcp}&& A' \doteq  \{R_{[t]}\leq w\}, \qquad B' \doteq    \{O_{x+y}\leq v\}, \qquad C' \doteq   \{Q_{[t],y}> z\},
\end{eqnarray}
where $w,v\in[0,\infty]$ and $z \in [-\infty,\infty)$ 
are arbitrary (the statistics $Y(t), Z(x+y), M(t,y)$
correspond to the L\'evy process $X$
in~\eqref{eq:Emebdding}, see~\eqref{eq:Levy_Overshoot}).
Since $X$
satisfies As.~\ref{A2} (as $S$ 
satisfies As.~\ref{A1}), the asymptotic independence in Thm.~ \ref{thm2} will follow from Thm.~\ref{thm} once we show that if $t$, $x$, and $y$ 
tend to infinity, in such a way that $t\wedge x\wedge y$ tends to infinity,\footnote{Here and throughout we use the notation: $a\vee b\doteq\max\{a,b\}$,
$a\wedge b\doteq\min\{a,b\}$ for any
$a,b\in\mbb R$.}
we have:
\begin{eqnarray}
\label{a}&& |P(A') - P(A)| \vee| P(B') - P(B)| \vee |P(C') - P(C)| = o(1),\\
\label{b}&& |P(A\cap B\cap C) - P(A'\cap B'\cap C')| = o(1).
\end{eqnarray}
Indeed, the triangle inequality implies
\begin{eqnarray*}
|P(A'\cap B'\cap C') - P(A')P(B')P(C')| &\leq& |P(A'\cap B'\cap C') - P(A\cap B\cap C)|\\ 
&+& |P(A\cap B\cap C) - P(A)P(B)P(C)|\\ &+& |P(A)P(B)P(C) - P(A')P(B')P(C')|,  
\end{eqnarray*}
which tends to zero if $t\wedge x\wedge y\to\infty$ in view of Eqns.~\eqref{a} 
and~\eqref{b} and  Thm.~~\ref{thm}\eqref{enum:i} (recall that $P(A\cap B\cap C) - P(A)P(B)P(C)=o(1)$ 
as $t\wedge x\wedge y\to\infty$).

\subsubsection{Proof of convergence in Eqn.~\eqref{a}} We proceed by showing 
the convergence 
of the three differences of probabilities by treating the three 
cases separately.  
Convergence of the first difference to zero 
follows by a duality argument. 
Since $S^*_n \doteq \max_{k\in\{0,\ldots,n\}} S_k$ 
increases to the random variable $S^*_\infty$ 
as $n\to\infty$ and $N(t)$ tends to infinity as $t\to\infty$
(both $P$-a.s.),
we have that the running supremum
$X^*(t)$ of $X$ also increases to $S^*_\infty$
as $t\to\infty$ $P$-a.s.  As a consequence, the indicators 
$\I_{\{X^*(t)\leq w\}}$ and $\I_{\{S^*_{[t]}\leq w\}}$ decrease to $\I_{\{S^*_{\infty}\leq w\}}$ as $t\to\infty$
$P$-a.s., 
and we have that both $P(X^*(t)\leq w)$ and $P(S^*_n\leq w)$ 
converge to $P(S^*_{\infty}\leq w)$. By duality, the random variables
$R_{[t]}$ and $Y(t)$ have the same laws as $S^*_{[t]}$ and $X^*(t)$, 
respectively, so 
that we find 
\begin{eqnarray}\label{eq:deltaA}
\Delta^A(t) \doteq P(A') - P(A)=o(1)\q\text{ as $t\to\infty$}.
\end{eqnarray}

The second difference is zero since definition~\eqref{eq:Emebdding}  yields 
$O_{x+y}=Z(x+y)$ and hence  the events 
$B$
and 
$B'$
coincide, implying the equality $P(B)=P(B')$ 
for all
$x,y>0$
and
$v\in[0,\infty]$.

In order to establish that the third difference in Eqn.~\eqref{a},
$\Delta^C(t,y)\doteq |P(C)-P(C')|$, tends to zero (for any
$z\in[-\infty,\infty)$)
as $t\wedge y$ tends to infinity, 
we need to control the deviation of the Poisson random variable 
$N(t)$ away from its mean $t$ as $t\uparrow\infty$. 
Pick $\delta\in(0,1/4)$, define the events
\begin{equation}
\label{eq:A_delta}
A_\d(t) \doteq \{N(t(1-\d)) < [t] < N(t(1+\d)) \},
\end{equation}
and note that the law of large numbers for subordinators~\cite[p.~92]{Bertoin} implies 
$P(A_\d(t))\to 1$ as $t\to\infty$. 
Recalling that $Q_{[t],y}=R^*_{[t]} - y$ and $M(t,y) = Y^*(t) - y$, 
we have 
\begin{eqnarray*}
\Delta^C(t,y) = 
|P(R^*_{[t]}> y+z, A_\delta(t))-P(Y^*(t)> y+z, A_\delta(t))| + o(1)
\quad\text{as $t\wedge y\to\infty$.}
\end{eqnarray*}
Since 
$\{R^*_n\}_{n\in\NN}$
is a non-decreasing process
and
$R^*_{N(s)}=Y^*(s)$,
$s\in\RR_+$,
the following holds:
\begin{eqnarray}
\hspace{10mm}
P(Y^*(t(1-\delta))> y+z, A_\delta(t))  \leq   
P(R^*_{[t]}> y+ z, A_\delta(t))
 \leq  P(Y^*(t(1+\delta))> y+z, A_\delta(t)).
\label{eq:star}
\end{eqnarray}
Hence, 
as
$t\wedge y\to\infty$,
we find 
\begin{eqnarray}
\nonumber
0\leq
\Delta^C(t,y) & \leq &
P(Y^*(t(1-\delta)) \leq y+z<Y^*(t(1+\delta)), A_\delta(t))
+ o(1)\\
\nonumber
& \leq &
P(t(1-\delta) \leq \tau(y+z)<t(1+\delta))+o(1)\\ \nonumber
& \leq &
P(\WH L(t(1-\delta)) \leq \WH L(\tau(y+z))\leq\WH L(t(1+\delta)))+o(1).
\end{eqnarray}
Eqn.~\eqref{eq:Lemma_iv_ineq} in Theorem~\ref{thm}\eqref{enum:iii}
implies 
\begin{equation}
\label{eq:cool_Equality2}
0\leq
\limsup_{t\wedge y\to\infty}\Delta^C(t,y)  
\leq \frac{8}{\te{}}\delta\qquad \forall \delta\in(0,1/4).
\end{equation}
Therefore we have
$\limsup_{t\wedge y\to\infty}\Delta^C(t,y) = \lim_{t\wedge y\to\infty}\Delta^C(t,y) =0$
and 
Eqn.~\eqref{a} follows.
\exit

\subsubsection{Proof of convergence in Eqn.~\eqref{b}} 
By the strong law of large numbers, 
the definition of $A_\d(t)$  
in~\eqref{eq:A_delta}
and the fact $B=B'$, the difference
$\Delta^*(t,x,y) \doteq P(A\cap B\cap C) - P(A'\cap B'\cap C')$ 
satisfies 
\begin{eqnarray}
\nn \Delta^*(t,x,y)
&=& P(A\cap B\cap C\cap (A'\cap C')^c, A_\d(t)) - 
P(A'\cap B\cap C'\cap (A\cap C)^c, A_\d(t)) + o(1)\\
\nn &=& P(A\cap B\cap C\cap [((A')^c\cap C')\cup (C')^c], A_\d(t)) \\ \nn &-& 
P(A'\cap B\cap C'\cap [(A^c\cap C)\cup C^c], A_\d(t)) + o(1) \\
\label{delta*} &=& P(A\cap (A')^c \cap E_\d(t)) - P(A'\cap A^c \cap E_\d(t)) \\
&+& P(C\cap (C')^c \cap F_\d(t)) - P(C'\cap C^c \cap F'_\d(t)) + o(1), 
\q\text{as $t\wedge x\wedge y\to\infty$},
\nn
\end{eqnarray}
where $E_\d(t) \doteq B\cap C\cap C'\cap A_\d(t)$, 
$F_\d(t) \doteq A\cap B\cap A_\d(t)$, 
$F'_\d(t) \doteq A'\cap B\cap A_\d(t)$ 
and
$\delta\in(0,1/4)$
is arbitrary but fixed.
The second difference of
the two probabilities  on the right-hand side of~\eqref{delta*}
satisfies the following estimates 
by the monotonicity of $R^*$ and $Y^*$ 
and the definition of $A_\d(t)$  
in~\eqref{eq:A_delta}
(cf. 
Eqn.~\eqref{eq:star} above):
\begin{eqnarray*}
\le|P(C\cap (C')^c \cap F_\d(t)) - P(C'\cap C^c \cap F'_\d(t))\ri| 
&\leq& P(C\cap (C')^c,A_\d(t)) + P(C'\cap C^c,A_\d(t))\\
&\leq& 2 P\le(Y^*(t(1-\delta)) \leq y+z<Y^*(t(1+\delta))\ri) + o(1).
\end{eqnarray*}
Hence, by an argument analogous to the one in~\eqref{eq:cool_Equality2},
this difference tends to zero as $t\wedge x\wedge y\to\infty$. 

Eqn.~\eqref{b}, and hence the asymptotic independence in 
Thm.~~\ref{thm2}, will follow once we verify the equality
\begin{eqnarray}\label{AAE}\q\ \ \ \ 
|P(A\cap (A')^c \cap E_\d(t)) - P(A'\cap A^c \cap E_\d(t))| = o(1)
\quad\text{as $t\wedge x\wedge y\to\infty$.}
\end{eqnarray}
The proof of the limit in~\eqref{AAE} is more involved than the one above as, unlike 
$Y^*$
and
$R^*$,
the processes
$Y$
and
$R$
do not have monotone paths. 
The main tools in this proof 
are the weak limit in~\eqref{eq:deltaA} and the strong Markov property of
$Y$.
Furthermore, to establish Eqn.~\eqref{AAE} we will need to control the time between 
the epochs $t$ and $\WH L^{-1}(\WH L(\tau(a))$ 
as
$t\wedge a\to \infty$
(for 
$a=x+y$
and
$a=z+y$), where
$\WH L^{-1}$
is the right-continuous inverse of $\WH L$
(see Section~\ref{subsec:Levy}
and~\cite[Ch.~IV]{Bertoin} for
definition). 

Consider the events $\tilde D_\d(t,a)$, $\bar D_\d(t,a)$ and
$D_\d(t,a) \doteq \tilde D_\d(t,a)\cup\bar D_\d(t,a)$
for any 
$a>0$:
\begin{eqnarray}
\label{eq:Dd}
\tilde D_\d(t,a) \doteq \le\{\WH L^{-1}(\WH L(\tau(a)))< (1-\delta) t\ri\}, & \quad &
\bar D_\d(t,a) \doteq \le\{\WH L^{-1}(\WH L(\tau(a))) > (1+\d) t\ri\}.
\end{eqnarray}
Since 
$D_\d^c(t,a) \subseteq 
\le\{\WH L(t(1-\d)) \leq \WH L(\tau(a)) \leq \WH L(t(1+\d))\ri\}$,
Theorem~\ref{thm}\eqref{enum:iii}
implies 
\begin{eqnarray}
\label{eq:Dco1}
0\leq
\limsup_{t\wedge a\to\infty}
P(D_\d^c(t,a))
\leq \frac{8}{\te{}}\delta\qquad \forall \delta\in(0,1/4).
\end{eqnarray}
Hence, by an argument analogous to~\eqref{eq:cool_Equality2},
equality~\eqref{AAE} 
will follow if we prove 
\begin{eqnarray}\label{AAE1}
|P(A\cap (A')^c \cap E_\d(t) \cap \tilde D_\d(t,x+y)) - P(A'\cap A^c \cap E_\d(t)\cap \tilde D_\d(t,x+y))| & = & o(1),\\
|P(A\cap (A')^c \cap E_\d(t) \cap \bar D_\d(t,x+y)) - P(A'\cap A^c \cap E_\d(t)\cap \bar D_\d(t,x+y))| & = & o(1),
\label{AAE2}
\end{eqnarray}
as $t\wedge x\wedge y\to\infty$.
With this in mind, we denote by  
$(\F_t)_{t\geq0}$
the completed right-continuous filtration 
generated by the process
$X$ in~\eqref{eq:Emebdding}.
Note that
$N$
is $(\F_t)$-adapted and define
an $(\F_t)$-stopping time
$T'_{[t]}\doteq\inf\{s\ge 0: N(s)=[t]\}$ 
(i.e.
$N(T'_{[t]})=[t]$).
Note further that 
by~\eqref{eq:Emebdding}
we have
$Y(T'_{[t]})=R_{[t]}$,
$Y^*(T'_{[t]})=R^*_{[t]}$
and
$Z(x+y)=O_{x+y}$
and hence the events $A', B', C' $ 
in~\eqref{abcp} can be expressed as:
$$
A'=\{Y(T'_{[t]}) \leq w\}, \q B'=\{Z(x+y)\leq v \}, \q C'=\{Y^*(T'_{[t]})> y+z\} = \{\tau(y+z) < T'_{[t]}\}.
$$
For any 
time
$t\geq0$
and
event $K$, 
which may depend on 
$t$,
define 
\begin{equation}
\label{eq:p_of_t}
p_K(t)\doteq P(K \cap A_\delta(t))
\end{equation}
and note that by~\eqref{eq:deltaA}
we have $\Delta^A(t)=p_{A\cap (A')^c}(t)-p_{A^c\cap A'}(t)+o(1)$
as $t\uparrow\infty$.
The key identity in our proof
is 
(for any
$a>0$)
given by
\begin{eqnarray}
\label{eq:p1_id}
\I_{\tilde D_\d(t,a)}P\left(\bar A\cap A_\delta(t)\left|\F_{\WH L^{-1}\le(\WH L(\tau(a))\ri)}\right. \right) & = &
\I_{\tilde D_\d(t,a)} p_{\bar A}\left(t-\WH L^{-1}\le(\WH L(\tau(a))\ri)\right)\qquad\text{$P$-a.s.} 
\label{eq:p2_id}
\end{eqnarray}
where
$\bar A$
denotes either
$A\cap (A')^c$ 
or
$A^c\cap A'$.
It is important to observe that, 
since the definition of 
$\bar A\cap A_\d(t)$
depends on the epoch
$t$
(cf.~\eqref{eq:A_delta} and the definitions of
$A, A'$ in~\eqref{abcp}),
for any
$\omega\in \tilde D_\d(t,a)$,
this set
on the right- (resp. left-) hand side of~\eqref{eq:p1_id}
is defined for the epoch
$t- \WH L^{-1}\le(\WH L(\tau(a))\ri)(\omega)$
(resp. $t$).
A formal construction of the right-hand side of~\eqref{eq:p1_id},
which ensures its measurability,  
requires the \textit{shift operator} $\theta$ that can be defined on the canonical probability
space and for any $s\in\mbb R_+$ 
and path $\w$ satisfies $\theta_{s}(\w)(\cdot) = \w(s + \cdot)$
(for details we refer to~\cite[Ch. 0]{Bertoin} and the references there in).  
The equality in~\eqref{eq:p1_id}
holds since 
$Y$
is a strong $(\F_t)$-Markov process, 
$Y\le(\WH L^{-1}\le(\WH L(\tau(a))\ri)\ri)=0$ 
and $A_\delta(t)$ is given by~\eqref{eq:A_delta}.
Note that the stopping time 
$\WH L^{-1}\le(\WH L(\tau(a))\ri)$
is the time when the excursion of 
$Y$,
during which
$Y$
crosses the level 
$a$
for the first time, ends.

\subsubsection{Proof of~\eqref{AAE1}}
\label{subsubsec:1}
We may assume
$x>z$
implying
$
\tau(z+y)\leq \tau(x+y)\leq  \WH L^{-1}\le(\WH L(\tau(y+x))\ri).
$
Hence
$B, \bar C, \tilde D_\d(t,x+y)\in\F_{\WH L^{-1}\le(\WH L(\tau(y+x))\ri)}$, 
where
$\bar C\doteq C\cap C' = \{\tau(y+z) < t\wedge T'_{[t]}\}$,
and~\eqref{eq:p1_id} yields
\begin{eqnarray}
\nn P(\bar A\cap B \cap \bar C \cap \tilde D_\d(t,y+x),A_\d(t)) 
&=& E\le[\I_{B\cap \bar C\cap\tilde D_\d(t,y+x)}P\left(\bar A\cap A_\delta(t)\left|\F_{\WH L^{-1}\le(\WH L(\tau(y+x))\ri)}\right. \right)\ri]\\
&=& E\le[\I_{B\cap \bar C\cap\tilde D_\d(t,y+x)}p_{\bar A}\left( t-\WH L^{-1}\le(\WH L(\tau(x+y))\ri)\right)\ri].
\label{eq:key_property_1}
\end{eqnarray}
The left-hand side of~\eqref{AAE1}
is by~\eqref{eq:key_property_1} bounded above by
\begin{eqnarray}
& & 
E\le[\I_{B\cap \bar C\cap\tilde D_\d(t,y+x)}\left\lvert 
p_{A\cap (A')^c}\left( t-\WH L^{-1}\le(\WH L(\tau(x+y))\ri)\right)
- p_{A^c\cap A'}\left( t-\WH L^{-1}\le(\WH L(\tau(x+y))\ri)\right)
\right\rvert
\ri].
\label{eq:key_property_2}
\end{eqnarray}
The expression in~\eqref{eq:key_property_2}
is 
$o(1)$
as $t\wedge x\wedge y\to \infty$ 
by~\eqref{eq:deltaA} (cf.~\eqref{eq:p_of_t})
and the dominated convergence theorem: for any sequence $(t_n, x_n, y_n)$,
such that $t_n\wedge x_n \wedge y_n\uparrow \infty$, 
pick an arbitrary subsequence $(t_{n_k}, x_{n_k}, y_{n_k})$ 
and note that, 
on the event 
$\cap_k \tilde D_\d(t_{n_k},y_{n_k}+x_{n_k})$,
the diferrence $t_{n_k} - \WH L^{-1}(\WH L(\tau(y_{n_k}+x_{n_k}))$
is bounded below by $\d\cdot t_{n_k}$ 
and hence tends to infinity. 
On the complement of this event, the random variable under the expectation
in~\eqref{eq:key_property_2}
is clearly zero for all large 
$k$.
This establishes~\eqref{AAE1}.


\subsubsection{Proof of~\eqref{AAE2}}
The proof in this case is slightly more complex than the one in Section~\ref{subsubsec:1}
as, intuitively speaking, it requires splitting the events in~\eqref{AAE2} in such a way so that 
the excursions of $Y$
corresponding to 
$\WH L(\tau(x+y))$
and 
$\WH L(\tau(z+y))$
are distinct (put differently, 
$Y$
crosses
levels 
$z+y$
and
$x+y$
for the first time during distinct excursions excursions away from $0$; recall that
$x>z$).

In order to analyse~\eqref{AAE2}, note 
$t(1+\delta)\leq\WH L^{-1}\le(\WH L(t(1+\delta))\ri)$
and hence
$\F_{t(1+\delta)} \subset\F_{\WH L^{-1}\le(\WH L(t(1+\delta))\ri)}$.
Therefore 
$\bar C\in \F_{t}$
and
$\bar A\cap A_\delta(t), \bar D_\d(t,x+y)\in \F_{t(1+\delta)}$
imply
$\bar A\cap A_\delta(t),\bar C, \bar D_\d(t,x+y)
\in \F_{\WH L^{-1}\le(\WH L(t(1+\delta))\ri)}$,
where
$\bar C=C\cap C'$
and
$\bar A$
as in~\eqref{eq:p1_id}.
The strong Markov property of 
$Y$
at the stopping time
$\WH L^{-1}\le(\WH L(t(1+\delta))\ri)$,
the equality
$Y\le(\WH L^{-1}\le(\WH L(t(1+\delta))\ri)\ri)=0$
and the definition of 
$B$
(cf.~\eqref{abcp})
yield
$$
P\le(B \cap H \left|\F_{\WH L^{-1}\le(\WH L(t(1+\delta))\ri)}\right.\ri)
=
\I_{H} P(B),\quad \text{
where
$H\doteq\{\WH L^{-1}(\WH L(\tau(y+x))) > \WH L^{-1}(\WH L(t(1+\delta)))\}$.}
$$ 
It holds
$\bar D_\d(t,x+y)\cap H^c\subseteq \{\WH L^{-1}(\WH L(\tau(y+x))) = \WH L^{-1}(\WH L(t(1+\delta)))\}$
and hence by Theorem~\ref{thm}\eqref{enum:iii}
we have 
$P(\bar D_\d(t,x+y)\cap H^c)=o(1)$ as $t\wedge x\wedge y\to \infty$. 
Therefore conditioning on 
$\F_{\WH L^{-1}\le(\WH L(t(1+\delta))\ri)}$
yields
\begin{eqnarray}
\nn
 P(\bar A\cap A_\d(t)\cap B \cap \bar C \cap \bar D_\d(t,y+x)) 
&=& 
 P(\bar A\cap A_\d(t)\cap \bar C \cap \bar D_\d(t,y+x), B \cap H  ) + o(1)\\
& = & 
P(\bar A\cap A_\d(t) \cap \bar C \cap \bar D_\d(t,y+x)) P(B)+o(1)
\label{eq:Prob_Important_not_final}
\end{eqnarray}
as $t\wedge x\wedge y\to \infty$.
Note that~\eqref{eq:Dd}, 
Theorem~\ref{thm}\eqref{enum:iii}
and the inclusion $\bar D_\d(t,y+z)\cap \bar C\subset \{\WH L(t)=\WH L(\tau(y+z))\}$ 
imply
$
P(\bar D_\d(t,y+z)\cap \bar C)\leq P(\WH L(t)=\WH L(\tau(y+z)))=o(1)
$
as $t\wedge y\to \infty$.
Hence, 
as $t\wedge x \wedge y\to \infty$, 
we can decompose the right-hand side of~\eqref{eq:Prob_Important_not_final} further 
using~\eqref{eq:p1_id} as follows:
\begin{eqnarray}\nonumber
\lefteqn{P(\bar A\cap A_\d(t) \cap \bar C \cap \bar D_\d(t,y+x))}\\
&=& 
\nonumber
P(\bar A\cap A_\d(t) \cap \bar C) -  P(\bar A\cap A_\d(t) \cap \bar C \cap \tilde D_\d(t,y+x))
- P(\bar A\cap A_\d(t) \cap \bar C \cap D^c_\d(t,y+x))\\
\nonumber
&=& 
P(\bar A\cap A_\d(t) \cap \bar C\cap \tilde D_\d(t,y+z)) -  P(\bar A\cap A_\d(t) \cap \bar C \cap \tilde D_\d(t,y+x))\\
\nonumber
& - & 
P(\bar A\cap A_\d(t) \cap \bar C\cap D^c_\d(t,y+z)) - P(\bar A\cap A_\d(t) \cap \bar C \cap D^c_\d(t,y+x)) + o(1)\\
 \label{aap2}
&=& E\le[\I_{\bar C\cap \tilde D_\d(t,y+z)}
p_{\bar A}\left( t-\WH L^{-1}\le(\WH L(\tau(z+y))\ri)\right)
- \I_{\bar C\cap \tilde D_\d(t,y+x)}p_{\bar A}\left( t-\WH L^{-1}\le(\WH L(\tau(x+y))\ri)\right)
\ri]\\
\nonumber
& - & 
P(\bar A\cap A_\d(t) \cap \bar C\cap D^c_\d(t,y+z)) - P(\bar A\cap A_\d(t) \cap \bar C \cap D^c_\d(t,y+x)) + o(1).
\end{eqnarray}
Since in~\eqref{aap2}
$\bar A$
denotes either
$A\cap (A')^c$ 
or
$A^c\cap A'$, 
the left-hand side in~\eqref{AAE2}
is 
bounded above by the sum
(recall $\Delta^A(t) = P(A\cap (A')^c) - P(A'\cap A^c)$, 
definition~\eqref{eq:p_of_t} and equalities~\eqref{eq:Prob_Important_not_final} and~\eqref{aap2}):
\begin{eqnarray*}
E\le[
\I_{\bar C\cap \tilde D_\d(t,y+z)}
\le|\Delta^A\left( t-\WH L^{-1}\le(\WH L(\tau(z+y))\ri)\right)\ri|+
\I_{\bar C\cap \tilde D_\d(t,y+x)}
\le|\Delta^A\left( t-\WH L^{-1}\le(\WH L(\tau(x+y))\ri)\right)\ri|
\ri] 
& & \\
+
2 P(D^c_\d(t,y+z)) + 2P(D^c_\d(t,y+x)) + o(1) 
\qquad \text{as $t\wedge x \wedge y\to \infty$.} 
& &  
\end{eqnarray*}
The expectation converges to zero 
as $t\wedge x \wedge y\to \infty$ 
by the same argument as in~\eqref{eq:key_property_2}
and, since 
$\delta\in(0,1/4)$ was arbitrary, the bound in~\eqref{eq:Dco1}
implies the equality in~\eqref{AAE2}.
\exit

\subsection{Proof of the Spitzer-type formula in~\eqref{eq:oo}}\label{ssec:spit} 
Let the compound Poisson process $X$ be as in~\eqref{eq:Emebdding}.
It is clear that the overshoots
$O_x$
of
$\{S_n\}_{n\in\NN^*}$
and
$Z(x)$
of
$X$
coincide
(see~\eqref{eq:Levy_Overshoot}
and Section~\ref{sec:Intro}
for definitions
of
$Z(x)$
and
$O_x$
respectively):
$Z(x)  =  O_x$
for any level
$x\in\RR_+.$
The ladder-height process of $X$ is
a compound Poisson process with Laplace exponent $\phi$
given by 
\begin{eqnarray*}
\frac{\phi(\th)}{\phi(0)} & = & \exp\le(\int_0^\infty \td t
\int_{[0,\infty)} (1 - \te{-\theta x})t^{-1} P(S_{N(t)}\in\td x)\ri)
 = 
\exp\le\{\sum_{n=1}^\infty\frac{1}{n}(1 - E[\te{- \th\, S_n^+}])\ri\}
\end{eqnarray*}
for any $\theta\in\mbb R_+$
(see e.g. \cite[p.166]{Bertoin}).
The second equality is a consequence of
Fubini's theorem and the fact
$\int_0^\infty t^{-1}P(N(t)=k)\td t 
= \Gamma(k)/k! = 1/k$, for all
$k\in\mbb N$,
where
$\Gamma$
denotes the Gamma function. 
This concludes the proof
of~\eqref{eq:oo}.~\exit
\vspace{-0.2cm}
\bibliographystyle{plain}
\bibliography{cite}

\end{document}